\documentclass{IEEEtran4PSCC}

\ifCLASSINFOpdf
   \usepackage[pdftex]{graphicx}
\else
   \usepackage[dvips]{graphicx}
\fi
\usepackage[cmex10]{amsmath}
\usepackage{booktabs}
\usepackage{multirow}
\usepackage{etoolbox}
\usepackage[intoc]{nomencl}
\makenomenclature
\renewcommand\nomgroup[1]{%
  \item[\bfseries
  \ifstrequal{#1}{A}{Indices and Sets}{
  \ifstrequal{#1}{B}{Sets}{
  \ifstrequal{#1}{C}{Parameters}{
  \ifstrequal{#1}{D}{Variables}{}}}}
]}
\usepackage{hyperref}
\usepackage{csquotes}
\usepackage{color}

\usepackage{tabularx}
\usepackage{mathtools}
\usepackage{algorithmic}
\usepackage{algorithm}
\usepackage{amssymb}
\usepackage{amsmath,amssymb}
\usepackage{enumitem}
\usepackage{amsthm}
\newtheorem{theorem}{Theorem}
\usepackage{makecell}

\begin{document}
%

\title{Successive Fixing for Large-Scale SCUC Using First-Order Methods}

\author{
\IEEEauthorblockN{Jinxin Xiong, Yanting Huang, Yingxiao Wang, Linxin Yang, Jianghua Wu, Shunbo Lei, Akang Wang} \\
\IEEEauthorblockA{Shenzhen Research Institute of Big Data, China \\
The Chinese University of Hong Kong (Shenzhen), Shenzhen, China
}
}

\maketitle

\begin{abstract}

Security-Constrained Unit Commitment is a fundamental optimization problem in power systems operations. The primary computational bottleneck arises from the need to solve large-scale Linear Programming (LP) relaxations within branch-and-cut. Conventional simplex and barrier methods become computationally prohibitive at this scale due to their reliance on expensive matrix factorizations. While matrix-free first-order methods present a promising alternative, their tendency to converge to non-vertex solutions renders them incompatible with standard branch-and-cut procedures. 
To bridge this gap, we propose a successive fixing framework that leverages a customized GPU-accelerated first-order LP solver to guide a logic-driven variable-fixing strategy. 
Each iteration produces a reduced Mixed-Integer Linear Programming~(MILP) problem, which is subsequently tightened via presolving. This iterative cycle of relaxation, fixing, and presolving progressively reduces problem complexity, producing a highly tractable final MILP model. When evaluated on public benchmarks exceeding 13,000 buses, our approach achieves a tenfold speedup over state-of-the-art methods without compromising solution quality.

\end{abstract}


\thanksto{\noindent 
This work is supported by Natural Science Foundation of China (Grant No.
12301416) and Guangdong Basic and Applied Basic Research Foundation (Grant No. 2024A1515010306). Corresponding author: Akang Wang $<$wangakang@sribd.cn$>$.}

\vspace{-2em}
\section{Introduction}
Security-Constrained Unit Commitment (SCUC) is a fundamental but complex power system optimization problem. 
The problem is notoriously challenging due to its high-dimensional, combinatorial nature and numerous security constraints that must be enforced. 
When formulated as a Mixed-Integer Linear Programming~(MILP) problem, significant advancements have been made through strengthening techniques, such as deriving tighter relaxations and employing advanced cutting planes~\cite{panConvexHullsUnit2017, knueven2020mixed}. These methods effectively reduce the computational burden by yielding a more compact branch-and-bound tree.
Consequently, modern SCUC formulations are often sufficiently tight in a sense that their Linear Programming (LP) relaxations provide solutions very close to the final integer optimum~\cite{kempke2025developing}. This characteristic implies that the overall efficiency of solving the SCUC model is fundamentally governed by the performance of solving its LP relaxations. 

Traditional methods like simplex and interior-point, while standard for SCUC LP relaxations, struggle with large-scale systems due to their reliance on expensive, hard-to-parallelize matrix factorizations. 
In contrast, \textit{first-order methods}~(FOMs)—such as the primal-dual hybrid gradient~\cite{lu2025cupdlp} and Halpern Peaceman-Rachford (HPR) methods~\cite{chen2025hpr}—replace these factorizations with efficient, highly parallelizable matrix-vector multiplications. 
This matrix-free design is ideal for GPU acceleration, positioning FOMs as a promising path toward scalable LP solutions~\cite{lu2025overview}.

However, a critical drawback is that FOMs typically yield \textit{non-vertex} LP solutions, which are not directly usable by MILP solvers. Core algorithms, such as branch-and-cut, require \textit{basic feasible solutions}~(i.e., vertices). 
Such solutions provide a basis structure that is essential for efficiently warm-starting the dual simplex method as well as deriving strengthening inequalities from the simplex tableau.
The work of~\cite{liu2024new} introduced crossover techniques to map non-vertex FOM solutions to vertices. However, the computational cost of this process grows significantly with problem size, as it is dominated by repeatedly solving large least-squares problems.
Consequently, for large-scale instances, the crossover phase alone can take longer than the entire initial process of finding a near-optimal solution.

An alternative strategy for leveraging FOM solutions in branch-and-cut involves their integration with established primal heuristics, such as feasibility pump~\cite{fischetti2005feasibility} and fix-and-propagate~\cite{mexi2023scylla}. A recent study by~\cite{kempke2025low} exemplifies this approach, applying FOM-derived solutions within a fix-and-propagate heuristic to the unit commitment problem without security constraints. Their framework first solves the full LP relaxation to inform the fixing of all binary variables, then optimizes the resulting LP for the continuous variables. However, this aggressive strategy of fixing all variables at once can lack robustness in finding high-quality feasible solutions. Moreover, their integration does not explicitly account for power system-specific structures or exploit the potential of modern GPU architectures. Finally, the application of FOMs to the more complex SCUC problem remains unexplored.
These critical gaps collectively motivate our central research question:

\begin{center}
    \textit{How can FOMs be effectively adapted to large-scale SCUC?}
\end{center}

\begin{figure*}[htbp]
\centering
\includegraphics[width=0.93\linewidth]{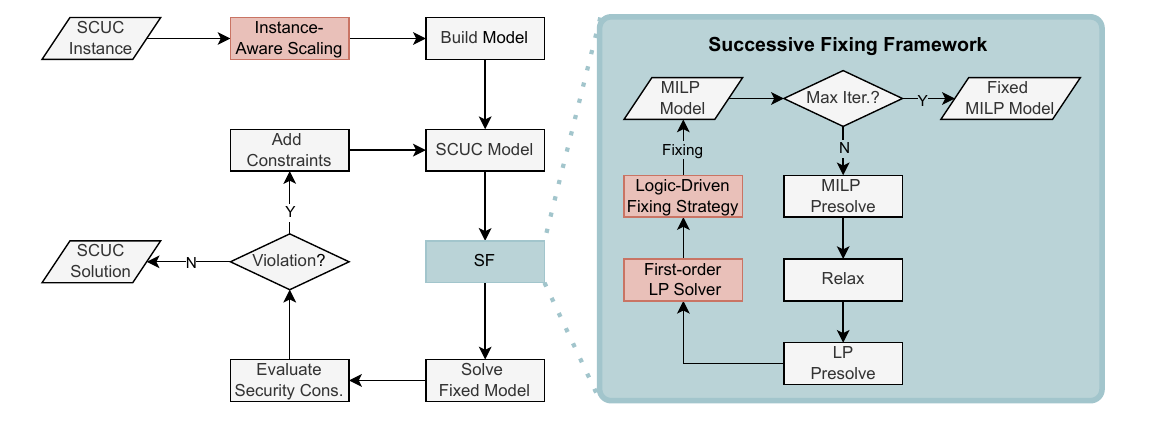}
\caption{A flowchart of our proposed framework, illustrating the transmission filtering outer loop (left panel) and the successive fixing inner loop (right panel). 
Key algorithmic enhancements are highlighted in red.} 
\label{fig:framework}
\end{figure*}

This work introduces a novel successive fixing framework to efficiently solve large-scale SCUC problems by leveraging FOMs. Our approach uses FOM solutions to guide a logic-driven variable-fixing strategy, progressively tightening and reducing the problem size through iterative presolving. 
This yields a final MILP model that is significantly more tractable. 
To enhance performance on SCUC LP relaxations, we introduce key algorithmic improvements that accelerate FOM convergence on modern GPU architectures without compromising solution quality. The complete framework is illustrated in Fig.~\ref{fig:framework}.
The distinct contributions of our work are as follows:

\begin{itemize}[leftmargin=*,nosep]
\item \textbf{Enhanced FOMs.} We introduce two algorithmic enhancements to the first-order LP solver HPR-LP—instance-aware preconditioning and low-precision GPU computation—that together substantially reduce the solution time for SCUC LP relaxations.

\item \textbf{Successive Fixing.} We propose a novel successive fixing framework for SCUC that integrates the HPR-LP solver with a logic-driven fixing strategy, enabling efficient and scalable solutions.

\item \textbf{Superior Performances.} 
Extensive experiments on public benchmarks, including systems with over 13,000 buses, show that our approach delivers high-quality solutions, achieving a $20\times$ speedup in LP relaxation and a $10\times$ acceleration in total solution time compared to baselines.
\end{itemize}

This paper is structured as follows. 
The SCUC problem is formulated in Section~\ref{sec:formulation}. 
Section~\ref{sec:FOM} then introduces the FOM solver and our key algorithmic improvements. 
We present the core successive fixing framework with logic-driven strategy in Section~\ref{sec:fixing}. 
Section~\ref{sec:experiment} evaluates the framework's performance through benchmarks and an ablation study. The paper concludes with a discussion of future work in Section~\ref{sec:conclusion}.

\section{Problem Formulation}
\label{sec:formulation}
\nomenclature[AB]{$g\in\mathcal{G}$}{Generator $g$ in the set of generators.}
\nomenclature[AB]{$t\in\mathcal{T}$}{Time index $t$ in the time periods $1,\cdots, T$.}
\nomenclature[AB]{$b\in\mathcal{B}$}{Bus $b$ in the set of buses.}
\nomenclature[AB]{$l\in\mathcal{L}$}{Line $l$ in the set of transmission lines.}
\nomenclature[AB]{$c\in\mathcal{C}$}{Contingency $c$ in the set of contingencies.}
\nomenclature[AB]{$\mathcal{G}_b$}{Set of generators connecting to bus $b$.}
\nomenclature[AB]{$h\in\mathcal{H}_g$}{Piecewise segment $h$ in the set of segments for $g$.}

\nomenclature[C]{$UT_g/DT_g$}{Minimum up/down time of generator $g$.}
\nomenclature[C]{$P^L_g/P^U_g$}{Minimum/maximum power limit of generator $g$.}
\nomenclature[C]{$P^U_{gh}$}{Maximum power for segment $h$ for generator $g$.}
\nomenclature[C]{$SU_g/SU_g$}{Start-up/down capacity of generator $g$.}
\nomenclature[C]{$RU_g/RD_g$}{Maximum ramp up/down capacity of generator $g$.}
\nomenclature[C]{$F_{l}^c$}{Thermal limit of line $l$ under 
contingency $c$.}
\nomenclature[C]{$C^{\text{pen}}$}{Penalty related to demand mismatch.}
\nomenclature[C]{$D_{bt}$}{Variable for load for bus $b$ in time $t$.}
\nomenclature[C]{$C_{gh}$}{Cost coefficient of segment $h$ of generator $g$.}
\nomenclature[C]{$C^L_{g}$}{Cost of generator $g$ operating at $P^L_g$.}

\nomenclature[D]{$u_{gt}$}{Commitment status of generator $g$ in time $t$.}
\nomenclature[D]{$v_{gt}$}{Start-up status of generator $g$ in time $t$.}
\nomenclature[D]{$w_{gt}$}{Shut-down status of generator $g$ in time $t$.}
\nomenclature[D]{$p'_{gt}$}{Power output above minimum by $g$ in time $t$.}
\nomenclature[D]{$p_{gth}$}{Power from segment $h$ for generator $g$ in time $t$.}
\nomenclature[D]{$\delta_{lb}^c$}{Power Transfer Distribution Factor (PTDF) of bus $b$ to line $l$ under contingency $c$.}
\nomenclature[D]{$d_{bt}^{\text{pen}}$}{Amount of unsatisfied demand for bus $b$ in time $t$.}
\nomenclature[D]{$c_{gt}$}{Production cost over $P^L_g$ for generator $g$ in time $t$.}
\printnomenclature

This section summarizes the key SCUC constraints and objectives. The complete mathematical model, including details on reserves, storage, and startup costs, can be found in~\cite{knueven2020mixed}.

\noindent\textit{Logic constraints:}
\begin{equation}
\label{eq:status}
\begin{aligned}
    & u_{gt}-u_{g,t-1}=v_{gt}-w_{gt} && \forall g\in\mathcal{G}, t\in \mathcal{T} \\
    & v_{gt} + w_{gt} \leq 1 && \forall g\in\mathcal{G}, t\in \mathcal{T}
\end{aligned}
\end{equation}

\noindent\textit{Minimum up/down:}
\begin{equation}
\label{eq:minup/down}
    \begin{aligned}
    & \sum_{\tau = \max{\{1, t - UT_g+1 \}}}^{t} v_{g\tau} \leq u_{gt} && \forall g \in \mathcal{G}, t \in \mathcal{T} \\
    & \sum_{\tau= \max{\{1, t - DT_g + 1\}}}^{t} w_{g\tau} \leq 1 - u_{gt} && \forall g \in \mathcal{G}, t \in \mathcal{T} 
\end{aligned}
\end{equation}

\noindent\textit{Production limits:}
\begin{equation}
    \label{eq:power_limit}
    p'_{gt} \leq (P^U_g - P^L_g) u_{gt} \;\; \forall g\in\mathcal{G}, t\in\mathcal{T}
\end{equation}

\noindent\textit{Ramp up/down:}
\begin{equation}
\resizebox{0.9\linewidth}{!}{$
    \label{eq:ramp up/down}
    \begin{aligned}
    & \left(p'_{gt} + P^L_g u_{gt}\right) - \left(p'_{g,t-1} + P^L_g u_{g,t-1}\right)\leq RU_gu_{g,t-1} + SU_gv_{gt} \\
    & \left(p'_{g,t-1}+ P^L_g u_{g,t-1}\right) - \left(p'_{gt} + P^L_g u_{gt}\right) \leq RD_gu_{gt} + SD_g w_{gt}  \\
     & \hspace{7cm} \forall g\in\mathcal{G}, t\in\mathcal{T} \\
    \end{aligned}
$}
\end{equation}

\noindent\textit{Production cost:}
\begin{equation}
    \label{eq:prod_cost}
    \resizebox{0.85\linewidth}{!}{$
    \begin{aligned}
    & p_{gth} \leq \left(P^U_{gh} - P^U_{g,h-1}\right) u_{gt} &&\forall g\in\mathcal{G}, t\in\mathcal{T}, h\in\mathcal{H}_g \\
    &\sum_{h\in\mathcal{H}_g}p_{gth} = p'_{gt} && \forall g\in\mathcal{G}, t\in\mathcal{T}\\
    & \sum_{h\in\mathcal{H}_g} C_{gh} p_{gth} = c_{gt} && \forall g\in\mathcal{G}, t\in\mathcal{T}
    \end{aligned}
    $}
\end{equation}

\noindent\textit{System-wide constraints:}
\begin{equation}
\resizebox{0.85\linewidth}{!}{$
    \begin{aligned}
    & \sum_{g\in\mathcal{G}} p_{gt} = \sum_{b\in\mathcal{B}} \left(D_{bt} - d_{bt}^{\text{pen}} \right) \quad\quad\quad \forall t \in \mathcal{T}\\
    & -F_l^c \leq \sum_{b\in\mathcal{B}} \delta_{lb}^c \left(\sum_{g\in\mathcal{G}_b} p_{gt} - \left(D_{bt}-d_{bt}^{\text{pen}}\right)\right) \leq F_l^c\\
    & \hspace{4cm}\forall l\in\mathcal{L}, t\in\mathcal{T}, c\in\{0\}\cup\mathcal{C} 
\end{aligned}
$}
\label{eq:security}
\end{equation}
The objective for SCUC is formulated as:
\begin{equation}
\label{eq:obj}
\sum_{g\in\mathcal{G}}\sum_{t\in\mathcal{T}}\left(c_{gt} + C^L_{g} u_{gt}\right) + C^{\text{pen}}\sum_{b\in \mathcal{B}} \sum_{t\in\mathcal{T}} d_{bt}^{\text{pen}}
\end{equation}

\section{Enhanced HPR-LP for SCUC Relaxations}
\label{sec:FOM}
When solving large-scale SCUC models via branch-and-bound, empirical evidence shows that optimizing the root node LP relaxation often becomes the computational bottleneck, to the point that it can take longer than the entire subsequent branch-and-bound search. This bottleneck motivates the adoption of recently developed FOMs, which leverage GPU parallelization to significantly accelerate the solution of large-scale LP problems~\cite{lu2025overview}.

This work employs the HRP method with semi-proximal terms for LP (denoted as \enquote{HPR-LP})~\cite{chen2025hpr} to solve large-scale SCUC relaxations. We consider LPs of the form:
\begin{equation} 
\label{eq:lp} 
	\min_{x \in \mathcal{K}} \; \left\{ \mu^\top x : A x \geq \theta \right\},
\end{equation}
where $\mathcal{K} \coloneqq \{ x \in \mathbb{R}^n \mid x^L \leq x \leq x^U \}$. Here, $x \in \mathbb{R}^n$ denotes the decision variables, $\mu \in \mathbb{R}^n$ the cost vector, and $x^L, x^U \in (\mathbb{R} \cup \{\pm\infty\})^n$ the variable bounds. The constraints are given by matrix $A \in \mathbb{R}^{m \times n}$ and right-hand side vector $\theta \in \mathbb{R}^m$, with $y \in \mathbb{R}^m$ denoting the associated dual variables.

The HPR method is an iterative first-order algorithm. Let $z$ be an auxiliary variable, $\Pi_{\mathcal{K}}(\cdot)$ the projection onto set $\mathcal{K}$, $\mathbb{I}(\cdot)$ the identity operator, $\sigma > 0$ a penalty parameter, and $\lambda \coloneqq \lambda_{\max}(AA^\top)$ a preconditioning parameter. Each iteration $k$ of HPR-LP consists of the following steps:
\vspace{-0.3em}
\begin{equation} 
\label{eq:hpr}
\resizebox{0.98\linewidth}{!}{$
\begin{aligned}
    &\bar{x}^{k+1} \leftarrow \Pi_\mathcal{K}\left(x^k+\sigma\left(A^\top y^k-\mu\right)\right) \\ 
    &\bar{y}^{k+1} \leftarrow \Pi_{\mathbb{R}_+^{m}}\left(y^k+\frac{1}{\lambda\sigma}\left(\theta-A\left(2 \bar{x}^{k+1}-x^k\right)\right)\right) \\ 
    &\bar{z}^{k+1} \leftarrow \frac{1}{\sigma}\left(\Pi_\mathcal{K} - \mathbb{I}\right) \left( x^k+\sigma\left(A^\top y^k-\mu\right)\right) \\ 
    &\left(\hat{x}^{k+1}, \hat{y}^{k+1}, \hat{z}^{k+1}\right) \leftarrow 2\left(\bar{x}^{k+1}, \bar{y}^{k+1}, \bar{z}^{k+1}\right)-\left(x^k, y^k, z^k\right) \\ 
    &\left(x^{k+1}, y^{k+1}, z^{k+1}\right) \leftarrow \frac{1}{k+2}\left(x^0, y^0, z^0\right)+\frac{k+1}{k+2}\left(\hat{x}^{k+1}, \hat{y}^{k+1}, \hat{z}^{k+1}\right) 
\end{aligned}
$}
\end{equation} 

The core HPR-LP updates for primal, dual, and auxiliary variables are given by the first three equations of~\eqref{eq:hpr}. Convergence is accelerated through two subsequent steps: Peaceman-Rachford relaxation~\cite{eckstein1992douglas} and Halpern iteration~\cite{lieder2021convergence}. 
A key advantage is that all core computations--including the estimation of $\lambda$ via the power method~\cite{strang2012linear}--rely entirely on matrix-vector multiplications, avoiding expensive matrix factorizations. This structure preserves the sparsity of constraint matrix $A$ and enables efficient GPU parallelization. As formalized below, HPR achieves $\mathcal{O}(1/k)$ iteration complexity in terms of the KKT residual.

\begin{theorem}[\cite{chen2025hpr}]
\label{th:conv}
    Assume the KKT solution for Problem~(\ref{eq:lp}) exists.
    Let $\{(\bar{x}^k, \bar{y}^k, \bar{z}^k)\}$ be the sequence generated by~(\ref{eq:hpr}), and $(x^*, y^*, z^*)$ be the corresponding limit point. 
    Then for all $k \geq 0$,
    \begin{equation*} 
\resizebox{0.98\linewidth}{!}{$
    \|\mathcal{R}\left((\bar{x}^{k+1}, \bar{y}^{k+1}, \bar{z}^{k+1})\right)\| \leq \frac{ R_0\left(\sigma\left(\|A\|+\left\|\sqrt{\lambda I_m - AA^\top} \right\|\right) + 1\right)}{\sqrt{\sigma}(k+1)},
    $}
    \end{equation*}
    where $
        \mathcal{R}(x,y,z) \coloneqq \left( 
        \begin{aligned}
            &y - \Pi_{\mathbb{R}_+^m}(y-Ax+\theta) \\ 
            & x - \Pi_\mathcal{K}(x-z) \\
            &\mu - A^\top y - z
        \end{aligned}
        \right)
        $
    denotes the KKT residual mapping for Problem~(\ref{eq:lp}), $I_m$ is the $m\times m$ identity matrix, and $R_0$ is a constant measuring the initial point's distance to the limit point.
\end{theorem}

To further enhance convergence, the full HPR-LP implementation features adaptive $\sigma$ updates and restart schemes. 
As a result, HPR-LP typically attains high-accuracy solutions faster than comparable GPU-accelerated first-order LP solvers~\cite{chen2025relationships}. 
Further details on HPR-LP and the proof of Theorem~\ref{th:conv} are available in~\cite{chen2025hpr,sun2025accelerating}.

\subsection{Instance-Aware Scaling}
\label{sec:scaling}
Theorem~\ref{th:conv} suggests that the convergence of HPR-LP is highly sensitive to the conditioning of the constraint matrix $A$, highlighting the critical role of preconditioning. While the default HPR-LP implementation uses iterative \textit{Ruiz scaling} for matrix equilibration, this process incurs a non-trivial computational overhead for large-scale problems. To mitigate this cost, we introduce a preemptive scaling strategy that is applied directly to the problem parameters prior to solving.

Analysis of the model formulation (Section~\ref{sec:formulation}) indicates that the source of poor conditioning is localized. Constraints with large-magnitude coefficients—namely, production limits~(\ref{eq:power_limit}), ramping constraints~(\ref{eq:ramp up/down}), and production cost definitions~(\ref{eq:prod_cost})—are driven by maximum capacity and cost parameters. Conversely, a significant portion of the model, including logic constraints~(\ref{eq:status}), minimum up/down constraints~(\ref{eq:minup/down}) and system-wide constraints~(\ref{eq:security}), is inherently well-scaled with unit coefficients.
Therefore, we introduce a production scaling parameter, $\eta^P$, and a cost scaling parameter, $\eta^C$, to normalize the problem data. 
These are defined as the maximum production limit and the maximum generation cost, respectively:
\begin{equation*}
    \eta^P \coloneqq \max_{g\in\mathcal{G}} P^L_g, \quad \eta^C \coloneqq \max_{g\in\mathcal{G}}\max_{h\in\mathcal{H}_g}C_{gh}.
\end{equation*}
We then scale all production-related parameters by $\eta^P$:
\begin{equation*}
\resizebox{0.95\linewidth}{!}{$
\begin{aligned}
    &\widetilde{P}^U_{gh} \coloneqq \frac{P^U_{gh}}{\eta^P},
    &&\widetilde{D}_{bt} \coloneqq \frac{D_{bt}}{\eta^P}, &&\widetilde{F}_l^c \coloneqq \frac{F_l^c}{\eta^P}, \\
    &\widetilde{RU}_g \coloneqq \frac{RU_g}{\eta^P},
    &&\widetilde{RD}_g \coloneqq \frac{RD_g}{\eta^P},
    &&\widetilde{SU}_g \coloneqq  \frac{SU_g}{\eta^P},
    &\widetilde{SD}_g \coloneqq  \frac{SD_g}{\eta^P}.
\end{aligned}
$}
\end{equation*}
And all cost-related parameters are scaled as follows:
$$
\widetilde{C}_{gh} \coloneqq \frac{\eta^P}{\eta^C}C_{gh}, \quad \widetilde{C}^\text{pen} \coloneqq \frac{\eta^P}{\eta^C} C^\text{pen}.
$$

This preconditioning generates a mathematically equivalent, better-conditioned formulation of the SCUC problem by rescaling its units, thereby improving numerical properties for FOMs. The binary commitment variables remain identical to the original problem, and the original production variables (e.g., $p'_{gt}, d_{bt}^{\text{curtail}}$) and objective value can be recovered by multiplying the scaled solutions by $\eta^P$ and $\eta^C$, respectively.

When using the preconditioned model, for HPR-LP, the computationally expensive Ruiz scaling is disabled, while the efficient Pock-Chambolle scaling~\cite{pock2011diagonal} and normalization of the right-hand-side and objective vectors are retained to enhance conditioning with minimal overhead.

\subsection{Accelerating Computation via Low-Precision Arithmetic}
While FOMs already benefit from GPU parallelization, we achieve substantial additional acceleration by employing single-precision (FP32) arithmetic instead of conventional double-precision (FP64)~\cite{markidis2018nvidia}. This optimization leverages two key advantages of modern GPU architectures:

\begin{itemize}[leftmargin=*,nosep]
\item \textit{Higher Throughput}: Modern GPUs have significantly more cores dedicated to single-precision operations, yielding a much higher theoretical throughput.
\item \textit{Reduced Memory Footprint}: Single-precision arithmetic halves the memory storage and bandwidth requirements, thereby alleviating a critical performance bottleneck.
\end{itemize}
The modest numerical imprecision introduced by FP32 arithmetic is effectively compensated by our preconditioning and scaling techniques (Section~\ref{sec:scaling}), which maintain solver stability and convergence. This precision trade-off is particularly advantageous in our successive fixing framework: since the LP solutions guide binary variable fixing rather than requiring exact optimality, the slight precision loss has a negligible impact on final solution quality. Meanwhile, the substantial computational acceleration enables more frequent LP solves within the fixing loop, dramatically improving the heuristic's overall throughput.

\section{Successive Fixing Framework}
\label{sec:fixing}
Although FOMs efficiently yield moderately accurate solutions to the LP relaxation, these solutions are typically non-vertex points. This contrasts with vertex solutions produced by simplex methods, which are more readily exploited by MILP solvers. While crossover procedures exist to recover a vertex solution from an interior point~\cite{liu2024new}, they often involve solving linear systems or least-squares subproblems, requiring matrix factorizations. This reintroduces the computational bottleneck that FOMs were chosen to avoid for large-scale problems.

To bridge this gap without sacrificing computational efficiency, we propose a successive fixing heuristic. Our approach leverages the well-documented inherent tightness of the 3-binary SCUC formulation and its strengthened variants~\cite{morales2013tight,ostrowski2011tight}, where the LP relaxation solution is often nearly integral. Starting from the FOM-generated solution, we apply a simple round-and-fix procedure. This method iteratively solves the LP relaxation and fixes binary variables with high confidence—those with values sufficiently close to 0 or 1 (e.g., beyond a threshold 
$\tau$)—to their rounded values. By progressively reducing the problem size, this framework transforms the original large-scale MILP formulation into a sequence of smaller LP problems, culminating in a final MILP model that is tractable for standard solvers. This strategy effectively bypasses the need for expensive matrix factorizations while capitalizing on the favorable structure of the SCUC problem.

\subsection{Logic-Driven Fixing}
Simple round-and-fix heuristics often produce infeasible intermediate MILP models, primarily due to violations of inter-temporal constraints such as the unit status constraints in \eqref{eq:status}. To address this, we introduce a specialized fixing strategy (Algorithm~\ref{alg:fixing_strategy}) designed to maintain feasibility with respect to core SCUC constraints throughout the solution process.

Our strategy is built upon two key operators. The first is a \emph{Round} operator, which takes a relaxed solution value $s$ and a confidence threshold $\tau \in [0, 0.5)$, mapping it to an integer or an indeterminate state:
$$
\operatorname{Round}(s, \tau) \coloneqq \begin{cases}
    1, & \text{if } s \geq 1-\tau \\
    0, & \text{if } 0\leq s \leq \tau \\
    -1, & \text{otherwise}.
\end{cases}
$$
A result of $-1$ indicates insufficient confidence to perform rounding. The second is a $\operatorname{Fix}(\mathcal{M}, S, s)$ operator, which fixes variable $S$ to value $s$ in model $\mathcal{M}$. For clarity, the pseudocode uses vectorized notation for element-wise operations.

The central innovation for preventing infeasibility---which naive independent rounding would cause---is a logic-driven consistency check (Algorithm~\ref{alg:fixing_strategy}, lines~\ref{step:criteria_start}--\ref{step:criteria_end}). This procedure ensures that for any generator and time step, the associated binary variables are fixed only if their relaxed values are both confident (i.e., $\operatorname{Round}(s,\tau) \neq -1$) and mutually consistent with the inter-temporal logic of constraints \eqref{eq:status}. By actively enforcing these temporal and physical constraints across the scheduling horizon, our method substantially improves the feasibility of the resulting fixed model.
\begin{algorithm}[h]
\caption{Fixing-Strategy }\label{alg:fixing_strategy}
\begin{algorithmic}[1]
\STATE \textbf{Input:} Relaxation solutions $\left\{\hat{u}, \hat{v}, \hat{w}\right\}_{\mathcal{G,T}}$, MILP model $\mathcal{M}$ and rounding threshold $\tau$. 

\FOR{$g\in\mathcal{G}$}
\FOR{$t=1,\cdots, T$}
\STATE $\hat{u}_{gt},\hat{v}_{gt},\hat{w}_{gt} \gets \operatorname{Round}\left(\{\hat{u}_{gt},\hat{v}_{gt},\hat{w}_{gt}\}, \tau\right)$
\ENDFOR

\FOR{$t=2,\cdots, T$} \label{step:criteria_start}
\IF{any of $\hat{u}_{gt}, \hat{u}_{g,t-1}, \hat{v}_{gt}, \hat{w}_{gt}$ is $-1$}
    \STATE Break.
\ENDIF
    \IF{$\hat{u}_{gt} - \hat{u}_{g,t-1} = \hat{v}_{gt} - \hat{w}_{gt}$ \textbf{and } $\hat{v}_{gt} +\hat{w}_{gt} \leq 1$}
        \STATE $\operatorname{Fix}(\mathcal{M}, \{u_{gt}, v_{gt}, w_{gt}\}, \{\hat{u}_{gt},\hat{v}_{gt}, \hat{w}_{gt}\})$
    \ELSE
        \STATE Break.
    \ENDIF
\ENDFOR \label{step:criteria_end}
\ENDFOR 
\RETURN $\mathcal{M}$  
\end{algorithmic}
\end{algorithm}

\subsection{Successive Fixing}
A single fixing round is often insufficient due to the conservative nature of the feasibility-preserving checks, which may leave many variables unresolved. We therefore employ an iterative successive fixing framework (Algorithm~\ref{alg:multiround}) that progressively reduces the problem complexity over multiple rounds.
In each round $i$, the current MILP model $\mathcal{M}$ is first tightened using a MILP presolve routine. 
This presolved model is then relaxed to an LP, which is further simplified by LP presolving. 
The resulting LP is solved efficiently using a first-order LP solver to obtain a relaxed solution $\left\{\hat{u}, \hat{v}, \hat{w}\right\}_{\mathcal{G,T}}$ for binary variables. This solution is passed to the feasibility-aware fixing strategy, Algorithm~\ref{alg:fixing_strategy}, which fixes a subset of confident and consistent variables within the original model $\mathcal{M}$.
This iterative process of presolving, solving, and fixing progressively reduces the problem size. The final, significantly smaller, MILP model is returned for the definitive integer solve.

\begin{algorithm}[htbp]
\caption{Successive Fixing Framework}\label{alg:multiround}
\begin{algorithmic}[1]
\STATE \textbf{Input:} MILP model $\mathcal{M}$, number of fixing rounds $R$

\FOR{ $i=1,\cdots, R$}
    \STATE $\widehat{\mathcal{M}}^i \gets \operatorname{MILP-Presolve}(\mathcal{M})$
    \STATE $\widehat{\mathcal{M}}^i \gets \operatorname{Relax}(\widehat{\mathcal{M}}^i)$
    \STATE $\widehat{\mathcal{M}}^i \gets \operatorname{LP-Presolve}(\widehat{\mathcal{M}}^i)$

    \STATE $\{\hat{u}^i, \hat{v}^i, \hat{w}^i\}_{\mathcal{G,T}} \gets \operatorname{FOM-LP}(\widehat{\mathcal{M}}^i)$
    \STATE $\mathcal{M}\gets \operatorname{Fixing-Strategy}(\{\hat{u}^i, \hat{v}^i, \hat{w}^i\}_{\mathcal{G,T}}, \mathcal{M})$
\ENDFOR

\RETURN $\mathcal{M}$
\end{algorithmic}
\end{algorithm}

\section{CASE STUDY}
\label{sec:experiment}
This section presents a comprehensive performance evaluation of the proposed method (denoted as \texttt{SF}) through two key analyses:
(i) Comparative benchmarking against state-of-the-art baselines, and
(ii) Ablation studies quantifying the individual contributions of the proposed components.
For reproducibility, our implementation is publicly available at~\url{https://github.com/jx-xiong/FOM-SCUC.git}.

\noindent\textbf{Baselines:}
We compare \texttt{SF} against two established SCUC solution strategies:

\begin{itemize}[leftmargin=*,nosep]
    \item \texttt{TF}: A monolithic approach where the full SCUC problem is solved directly using Gurobi, augmented with the \textit{transmission filtering} technique from~\cite{xavier2019transmission} to handle security constraints.
    
    \item \texttt{TD}: A \textit{temporal decomposition} approach implemented in UnitCommitment.jl~\cite{xavier2022unitcommitment}, which partitions the problem into sequential subproblems solved using Gurobi with transmission filtering. The subproblem duration and advancement window are both set to 6 time intervals.
\end{itemize}

\noindent\textbf{Configuration.}
The proposed framework is implemented as an extension of UnitCommitment.jl~\cite{xavier2022unitcommitment}. All experiments were conducted in Julia v1.10.4, using Gurobi Optimizer v11.0.1~\cite{gurobi} as the underlying MILP solver and HPR-LP v0.1.0~\cite{chen2025hpr} as the first-order LP solver. 
Computations were performed on a server with a 13th Gen Intel(R) Core(TM) i9-13900K processor and NVIDIA RTX 4090 GPU, using 4 parallel threads for both transmission filtering and Gurobi. Each instance was subject to a 3,600-second time limit.

We employ the two-stage transmission filtering strategy from~\cite{xavier2022unitcommitment}, beginning with a $1\%$ optimality gap. If no security violations are detected, the gap is tightened to $0.1\%$. Correspondingly, our fixing framework parameters are staged: for the initial $1\%$ gap phase, we set confidence threshold $\tau=0.1$ with $R=2$ fixing rounds; for the refined $0.1\%$ gap phase, we increase to $R=4$ rounds.

\noindent \textbf{Benchmark.}
We evaluate the proposed framework using instances from the MATPOWER dataset~\cite{zimmerman2010matpower} with $T=36$ time periods. 
Our analysis focuses on the 20 largest instances in the dataset, each comprising over 1,000 buses.

\noindent \textbf{Performance Metrics.}
We employ the following metrics to assess solution quality and computational efficiency:
\begin{itemize}[leftmargin=*,nosep]
    \item Relative Gap (Rel. Gap): 
    The percentage difference between a method's objective value ($\nu$) and that of the \texttt{TF} baseline $\nu^{\texttt{TF}}$, calculated as $\frac{\nu -\nu^{\texttt{TF}}}{\nu^{\texttt{TF}}} \times 100\%.$

    \item Time Ratio: The solution time ($s$) relative to the \texttt{TF} baseline $s^{\texttt{TF}}$, given by the ratio $s/s^{\texttt{TF}}$.

    \item $\text{SGM}_{10}$: Scaled shifted (by 10 seconds) geometric mean of runtimes.
\end{itemize}

\subsection{Comparing Against Baselines}
In this section, we compare the performance of our proposed method, \texttt{SF}, against the \texttt{TF} and \texttt{TD} baselines. The results are summarized in Table~\ref{tab:result}. Column \enquote{Obj. (Rel. Gap)} reports the objective value (scaled by $10^7$) alongside its relative gap to the \texttt{TF} baseline. The \enquote{Time} column presents the total solution time and its ratio to the \texttt{TF} baseline. The bottom section of the table summarizes the number of instances solved to feasibility by each method, along with the average performance metrics computed over the subset of instances that all methods successfully solved.
\begin{table}[!ht]
\centering
\caption{Results comparing \texttt{SF} against \texttt{TF} and \texttt{TD}.}
\label{tab:result}
\renewcommand{\arraystretch}{1.2} 
    {\LARGE
\resizebox{\linewidth}{!}{
\begin{tabular}{lllllll}
\toprule
\multirow{2}{*}{\bf Instance}& \multicolumn{2}{c}{\texttt{TF}} &\multicolumn{2}{c}{\texttt{TD}}& \multicolumn{2}{c}{\texttt{SF}}  \\
        \cmidrule(l){2-3}\cmidrule(l){4-5}\cmidrule(l){6-7}
        & Obj. & Time (s) & Obj. (Rel. Gap) & Time & Obj. (Rel. Gap) & Time\\
\midrule
1354pegase & 1.575 & 18 & 1.614 (2.49\%) & 13 (0.72) & 1.575 (0.00\%) &    7 (0.39) \\
1888rte & 2.345 & 100 & 2.385 (1.70\%) & 29 (0.29) & 2.346 (0.01\%) &    19 (0.19) \\
1951rte & 2.494 & 136 & 2.531 (1.46\%) & 19 (0.14) & 2.494 (0.01\%) &   10 (0.07) \\
2383wp & 1.369 &   10 & 1.381 (0.92\%) & 16 (1.60) & 1.368 (0.00\%) & 14 (1.40) \\
2736sp & 0.962 & 15 & 0.994 (3.34\%) & 17 (1.13) & 0.966 (0.42\%) &   10 (0.67) \\
2737sop & 0.847 & 9 & 0.893 (5.48\%) & 16 (1.78) & 0.854 (0.83\%) &   5 (0.56) \\
2746wop & 0.853 & 51 & - & - & 0.857 (0.46\%) &   17 (0.33) \\
2746wp & 0.965 & 12 & 1.002 (3.85\%) & 16 (1.33) & 0.972 (0.74\%) &   7 (0.58) \\
2848rte & 2.427 & 477 & 2.463 (1.47\%) & 60 (0.13) & 2.430 (0.14\%) &   48 (0.10) \\
2868rte & 2.497 & 288 & 2.529 (1.26\%) & 39 (0.14) & 2.497 (0.01\%) &   16 (0.06) \\
2869pegase & 3.895 & 232 & 4.007 (2.89\%) &   69 (0.30) & 3.896 (0.02\%) & 75 (0.32) \\
3012wp & 1.204 & 39 & 1.243 (3.19\%) & 23 (0.59) & 1.206 (0.10\%) &   17 (0.44) \\
3120sp & 1.177 & 65 & 1.210 (2.85\%) & 33 (0.51) & 1.179 (0.17\%) &   26 (0.40) \\
3375wp & 1.418 & 672 & 1.544 (8.93\%) &   57 (0.08) & 1.419 (0.10\%) & 76 (0.11) \\
6468rte & 5.706 & 784 & 5.830 (2.18\%) &   293 (0.37) & 5.706 (0.01\%) & 342 (0.44) \\
6470rte & - & - & 6.744 &   525 & 6.643 & 1,572 \\
6495rte & 5.403 & 2,335 & 5.519 (2.14\%) &   252 (0.11) & 5.405 (0.04\%) & 516 (0.22) \\
6515rte & 5.573 & 754 & 5.697 (2.22\%) & 287 (0.38) & 5.574 (0.03\%) & 485 (0.64) \\
9241pegase & - & - & 11.599 & 1,441 & 11.357 &   1,118 \\
13659pegase & 26.561 & 1,776 & 27.004 (1.67\%) & 688 (0.39) & 26.564 (0.01\%) &  203 (0.11) \\
\midrule 
\textbf{Count} &18 & &19 & & \bf 20 & \\
\textbf{Avg.} && & \hspace{1.7cm}(2.82\%)&\hspace{1.cm} (0.59) & \hspace{1.5cm}(\textbf{0.15\%}) &\hspace{1.cm}(\textbf{0.39})\\
\bottomrule
\end{tabular}
}
}
\renewcommand{\arraystretch}{1} 
\end{table}

The monolithic \texttt{TF} method achieves the best objective value on nearly all instances and demonstrates competitive solution times for small-to-moderate problems, even being the fastest on instance \enquote{2383wp}. However, its performance degrades substantially with increasing problem scale, failing to identify feasible solutions for two of the largest instances within 3,600 seconds.
In contrast, the \texttt{TD} method exhibits superior scalability. It requires, on average, approximately 60\% of the time taken by \texttt{TF} and is the fastest method for several large instances. For example, on \enquote{6495rte}, it terminates and generates a SCUC solution within a tenth of the time required by \texttt{TF} and half the time of \texttt{SF}. This computational advantage, however, comes at the expense of solution quality. On instance \enquote{3375wp}, \texttt{TD} yields a solution with a 8.93\% relative gap despite being significantly faster. In a more severe case (\enquote{2746wop}), the method's myopic decomposition leads to an infeasible subproblem, preventing it from finding any feasible solutions.

The proposed \texttt{SF} method effectively synthesizes the strengths of both baselines, achieving robust performance in both solution quality and computational efficiency. A key advantage over the \texttt{TD} baseline is that \texttt{SF} successfully obtained feasible solutions for all 20 test instances. Although \texttt{SF} was slower than \texttt{TD} on some moderate-to-large systems (e.g., from \enquote{3375wp} to \enquote{6515rte}), it consistently delivered solution quality nearly on par with \texttt{TF}. This is evidenced by a minimal average optimality gap of only 0.15\% across all instances, far surpassing \texttt{TD}'s performance.

The efficiency of \texttt{SF} is particularly pronounced when compared to the \texttt{TF} baseline. On multiple instances, including \enquote{1951rte}, \enquote{2848rte}, and the large-scale \enquote{13659pegase}, our method achieved significant speedup factors ranging from $5\times$ to approximately $20\times$. For the two largest cases, \enquote{9241pegase} and \enquote{13659pegase}, \texttt{SF} demonstrated superior speed \emph{and} solution quality compared to the \texttt{TD} method. This combination of high efficiency and minimal optimality loss underscores the practical value of our approach.

\subsection{Ablation Study}
To quantify the impact of each component in our framework, we perform an ablation study. The analysis starts from a minimal implementation, with components added sequentially to construct the final \texttt{SF} method. 
The specifications for each ablated method are provided in Table~\ref{tab:ablation_method}.

\vspace{-1em}
\begin{table}[!ht]
    \centering
    \caption{Configuration of ablation methods.}
    \label{tab:ablation_method}
\renewcommand{\arraystretch}{1.05} 
    \resizebox{1\linewidth}{!}{
    \begin{tabular}{lcccc}
        \toprule
           \bf Method & \bf \makecell{ Successive\\Fixing} & \bf LP Solver & \bf \makecell{Instance\\Scaling} & \bf Precision\\
        \midrule
        \texttt{SF(Gurobi)} & Yes & Gurobi  & No & -\\
        \texttt{SF(FOM)} & Yes & HPR-LP  & No & FP64\\
        \texttt{SF(FOM+FP32)} & Yes & HPR-LP & No & FP32 \\
        \texttt{\texttt{F}} & No & HPR-LP & Yes & FP32\\
        \texttt{\texttt{SF}} & Yes & HPR-LP & Yes & FP32\\
        \bottomrule
    \end{tabular}
    }
\renewcommand{\arraystretch}{1.} 
\end{table}

\begin{figure}[!ht]
    \centering
    \includegraphics[width=1.\linewidth]{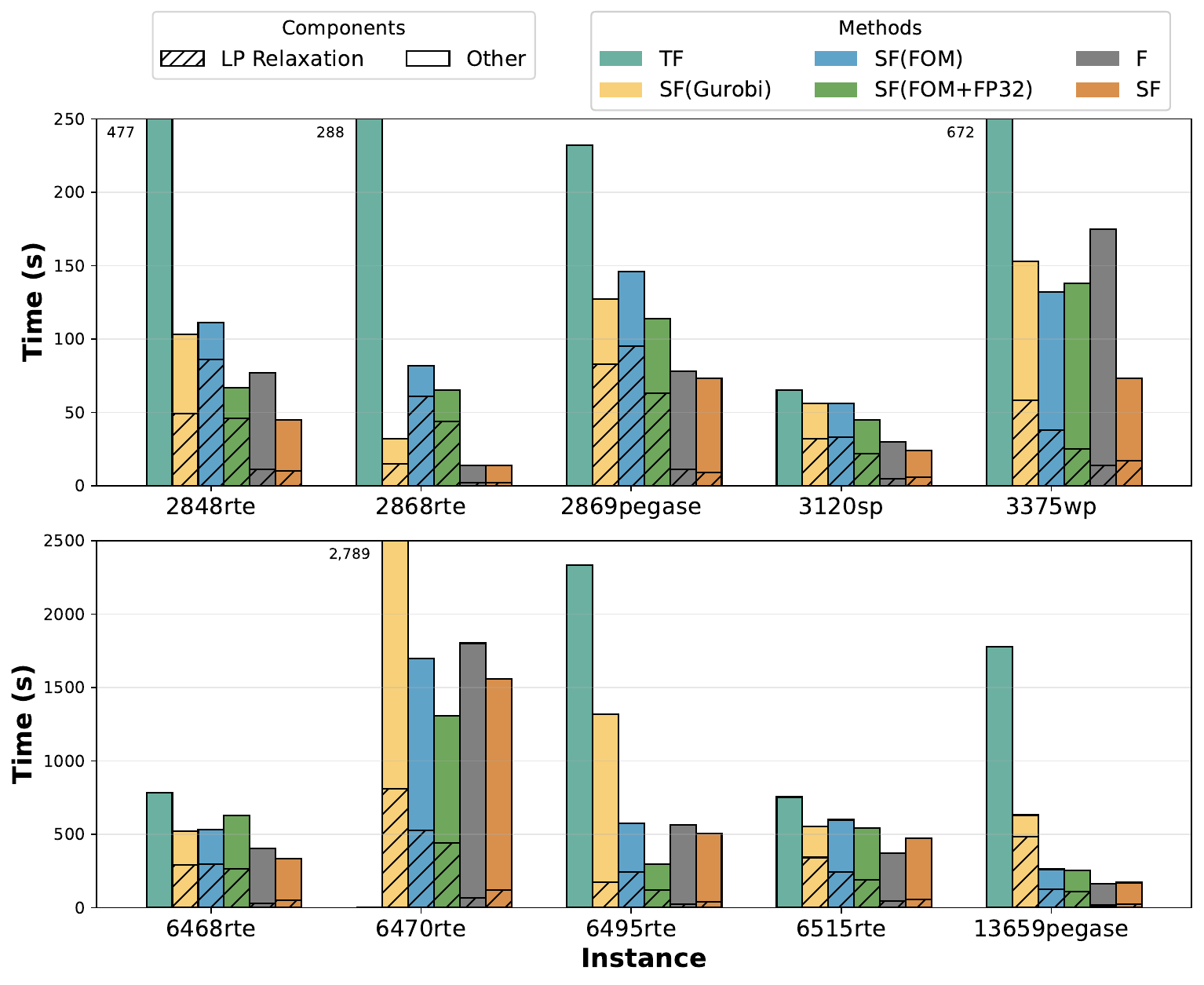}
    \vspace{-0.5cm}
    \caption{Time decomposition for the ablation study. Each bar represents the total solution time, containing the time for solving LP relaxations and other procedures (such as presolve, fixing, branching and filtering). }
    \label{fig:decompose}
\end{figure}

\begin{figure}[!ht]
    \centering
    \includegraphics[width=0.95\linewidth]{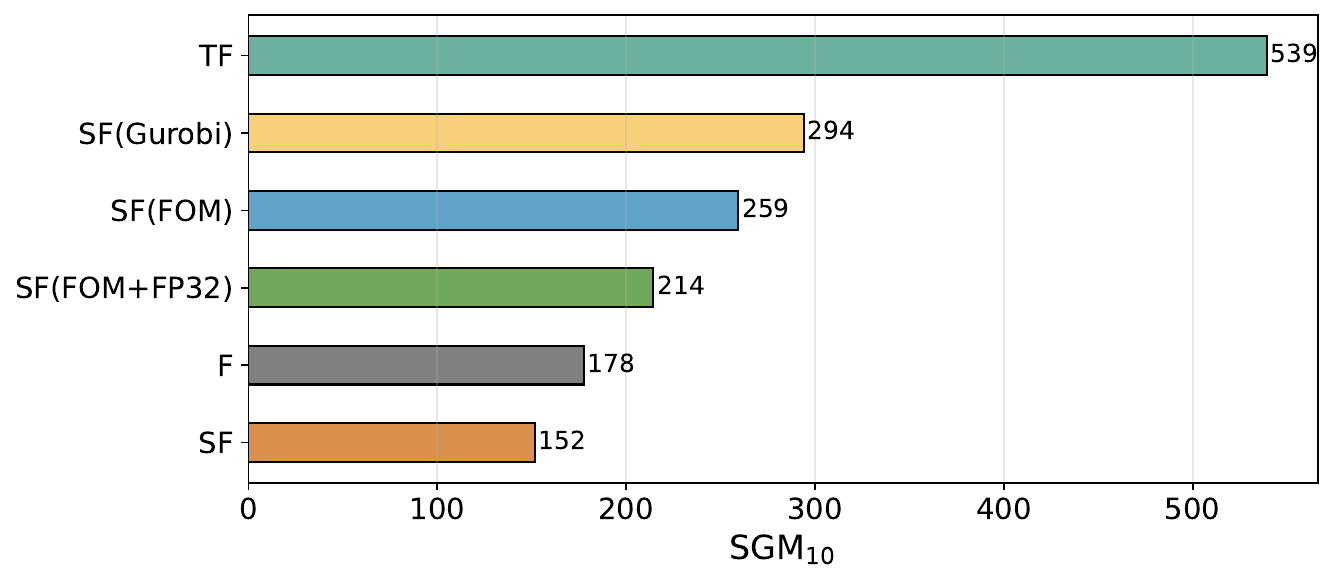}
    \vspace{-0.2cm}
    \caption{$\operatorname{SGM_{10}}$ of the \texttt{TF} and the ablation methods across the 10 largest instances solved by all ablation methods. (\enquote{6470rte} is excluded for \texttt{TF}.) }
    \label{fig:sgm_time}
\end{figure}

To ensure a fair performance comparison, this section focuses on the 10 largest instances (shown in Fig.~\ref{fig:decompose}) that were successfully solved by all ablation methods within the time limit. 
Fig.~\ref{fig:decompose} breaks down the total solution time, distinguishing between the cost of solving LP relaxations and other procedures (e.g., presolve, fixing, branching, and transmission filtering). 
Complementing this, Fig.~\ref{fig:sgm_time} reports the $\operatorname{SGM_{10}}$ runtime for this instance subset.


The results confirm that our successive fixing framework yields significant speedups over the monolithic \texttt{TF} baseline, regardless of the underlying LP solver. 
However, the choice and configuration of LP solvers are critical to performance. 
The naive integration of the first-order method (\texttt{SF(FOM)}) is slower than the Gurobi-based approach (\texttt{SF(Gurobi)}) on moderately scaled instances. This gap is likely due to overhead from default iterative scaling and the superior efficiency of state-of-the-art MILP solvers at this problem size.
The advantage of deploying FOMs becomes decisive at larger scales, as demonstrated by the cases of \enquote{6470rte}, \enquote{6495rte}, and \enquote{13659pegase}. On the \enquote{13659pegase} system, the enhanced HPR-LP solver achieves a $20\times$ speedup for LP relaxations, reducing the solution time from 485 seconds with \texttt{SF(Gurobi)} to just 22 seconds, while the time for other processes remains comparable.

Subsequent enhancements to the FOM integration further improve performance. Employing single-precision computation on the GPU for core FOM steps reduces the $\operatorname{SGM_{10}}$ runtime from 259 to 214. 
A customized scaling method provides an additional gain, lowering the metric to 152, a 41\% total reduction. 
The efficiency of these algorithmic enhancements is further illustrated in Fig.~\ref{fig:decompose}, where the components for LP relaxations for \texttt{F} and \texttt{SF} are consistently the shortest across all instances. 
On average, the complete \texttt{SF} method achieves a speedup of approximately 70\% compared to \texttt{TF} and 50\% compared to \texttt{SF(Gurobi)}.

Finally, a comparison between the single round fixing \texttt{F} and the full successive fixing framework (\texttt{SF}) highlights the latter's effectiveness. For example, on \enquote{3375wp}, \texttt{SF} incurs a slightly higher cost in solving LP relaxations than \texttt{F} (17 seconds versus 14 seconds). However, this is offset by a drastic reduction in time for solving the MILP problems (40 seconds versus 145 seconds). These results demonstrate that the initial investment in multi-round LP solutions, fixing, and presolving significantly reduces the complexity of the final MILP model, thereby validating the efficiency of the proposed \texttt{SF} framework.

\section{Conclusions}
\label{sec:conclusion}
This paper presented a novel successive fixing framework to tackle the computational bottleneck of large-scale SCUC. The core innovation lies in effectively leveraging non-vertex solutions from a GPU-accelerated FOM to guide a logic-driven fixing strategy for progressively simplifying the underlying MILP model. Key algorithmic enhancements—including instance-aware preconditioning and low-precision computation on GPUs—were critical to making the FOM solver practical and efficient for this role.
On large-scale benchmarks, our framework achieves significant speedups over state-of-the-art solvers while maintaining high solution quality. 
Future work will explore deeper integration of the FOM solution to better prune the MILP search space and improve high-precision performance.

\vspace{-0.24cm}
\bibliographystyle{IEEEtran}
\bibliography{ref}

\end{document}